\documentclass[12pt]{article}
\usepackage{mathrsfs}
\usepackage{amsmath,amsfonts,amssymb,rotating,amsthm}
\usepackage{hyperref}
\usepackage{array}
\usepackage{breqn}
\usepackage{color}
\usepackage{fancyheadings}
\usepackage{mathptmx}
\usepackage{dsfont}
\usepackage{caption}
\usepackage{subcaption}
\usepackage{tikz}
\usetikzlibrary{trees}

\tikzstyle{every node}=[circle,inner sep=1pt,fill=white!60]
\tikzstyle{tn}=[shape=circle, draw, color=black!70]
\tikzstyle{marke}=[shape=circle,minimum size=0.2cm, draw,blue]

\allowdisplaybreaks
\def\qed{\nopagebreak\hfill{\rule{4pt}{7pt}}}
\def\proof{\noindent {\it{Proof.} \hskip 2pt}}
\newtheorem{thm}{Theorem}[section]

\theoremstyle{remark}

\numberwithin{equation}{section}

\hfuzz=\maxdimen
\begin{document}

\begin{center}
{\large\bf A Context-free Grammar for the
$e$-Positivity of

the Trivariate Second-order
Eulerian Polynomials
}

\vskip 6mm

William Y.C. Chen$^1$ and Amy M. Fu$^2$

\vskip 3mm

$^{1}$Center for Applied Mathematics\\
Tianjin University\\
Tianjin 300072, P.R. China

\vskip 3mm

$^{2}$School of Mathematics\\
Shanghai University of Finance and Economics\\
Shanghai 200433, P.R. China

 Email: {\tt chenyc@tju.edu.cn, fu.mei@mail.shufe.edu.cn}

\end{center}

\noindent{\bf Abstract.} Ma-Ma-Yeh made a beautiful
observation that  a transformation of the grammar of Dumont
 instantly leads to the $\gamma$-positivity of the
Eulerian polynomials. We notice that the transformed
 grammar bears a striking resemblance to the grammar for 0-1-2 increasing trees also due to Dumont.
The appearance of the factor of two fits perfectly
 in a grammatical labeling of 0-1-2 increasing plane trees. Furthermore,
 the grammatical calculus
  is instrumental to the computation
of the generating functions. This approach can be adapted
to study the $e$-positivity of the trivariate
 second-order Eulerian
 polynomials first introduced by Dumont in the
 contexts of ternary trees and Stirling permutations, and
 independently defined by Janson, in connection with the
joint distribution of the numbers of ascents,
  descents and     plateaux over
Stirling permutations.

\vskip 3mm


\vskip 3mm


\section{Introduction}

The objective of this paper is to present a
context-free grammar to derive
the $e$-positivity of the trivariate
second-order Eulerian polynomials $C_n(x, y, z)$
defined
on  Stirling permutations, first introduced
by Dumont~\cite{Dumont-1980} in terms of ternary
trees and Stirling permutations, and rediscovered
by Janson~\cite{Janson-2008}.

This work was inspired by a beautiful
 observation of Ma-Ma-Yeh~\cite{Ma-Ma-Yeh-2019}
 that   a transformation of
 a context-free grammar found by Dumont~\cite{Dumont-1996}
instantly leads to the $\gamma$-positivity
of the Eulerian polynomials. We find that the transformed
grammar not only implies the $\gamma$-positivity, it
also provides a combinatorial interpretation of the
$\gamma$-coefficients in terms of
increasing plane trees.

For $n\geq 1$, let $[n]=\{1,2, \ldots, n\}$ and let
$S_n$ denote the set of permutations of $[n]$.
For a permutation $\sigma=\sigma_1\sigma_2 \cdots \sigma_n
\in S_n$, we assume that a zero is patched at the
beginning and at the end, that is, $\sigma_0=\sigma_{n+1}=0$. An index $1\leq i \leq n$
is said to be a descent (ascent) of a permutation $\sigma \in S_n$ if $\sigma_i > \sigma_{i+1}$  $(\sigma_{i-1} <
\sigma_{i})$. The number of permutations of $[n]$ with
$k$ descents $(1\leq k \leq n)$ is often denoted by
$A(n,k)$, or sometimes by $\left\langle {n \atop k}  \right\rangle$. The Eulerian polynomials
$A_n(x)$ are defined by $A_0(x)=1$ and for $n\geq 1$,
\begin{equation}
A_n(x) = \sum_{\sigma\in S_n} x ^{{\rm des } (\sigma)} = \sum_{k=1}^n A(n,k)x^k,
\end{equation}
where ${\rm des }(\sigma)$ denotes the number
of descents of a permutation $\sigma$. A bivariate version
of the Eulerian polynomials is given by
\begin{equation}
A_n(x,y) = \sum_{\sigma\in S_n} x ^{{\rm des } (\sigma)}
y ^{{\rm asc } (\sigma)} = \sum_{k=1}^n A(n,k)x^ky^{n+1-k},
\end{equation}
where $n\geq 1$ and ${\rm asc} (\sigma)$
 stands for the number of ascents of $\sigma$. Bear in mind that for any permutation
$\sigma \in S_n$, we have
\begin{equation} {\rm des } (\sigma) +
{\rm asc } (\sigma) = n+1.
\end{equation}

The first few
 values of $A_n(x,y)$ are given below,
\begin{eqnarray*}
A_1(x,y) & = & xy, \\[6pt]
A_2(x,y) & = & xy^2+x^2y, \\[6pt]
A_3(x,y) & = & xy^3+4x^2y^2+x^3y, \\[6pt]
A_4(x,y) & = & xy^4 + 11 x^2y^3 + 11 x^3 y^2 + x^4y,\\[6pt]
A_5(x,y) & = & xy^5+ 26x^2y^4 +66x^3y^3 +26 x^4y^2 + x^5y, \\[6pt]
A_6(x,y) & = & xy^6+57x^2y^5+302x^3y^4+ 302 x^4y^3 + 57 x^5y^2 + x^6y.
\end{eqnarray*}

A celebrated theorem
of Foata and Sch\"utzenberger~\cite{Foata-Schutzenberger-1970} states that for $n
\geq 1$, the Eulerian polynomial $A_n(x)$ can be
expanded uniquely in the following form
\begin{equation} \label{A-n-gamma}
   A_n(x) = \sum_{k=1}^{[(n+1)/2]}
    \gamma_{\,n,\,k}  \,  x^k(1+x)^{n-2k+1}
\end{equation}
with nonnegative coefficients $\gamma_{\,n,\,k}$.
The above expression (\ref{A-n-gamma}) is
called the $\gamma$-expansion of $A_n(x)$, which
can be restated  as
\begin{equation} \label{gamma-x-y}
A_n(x,y) = \sum_{k=1}^{[(n+1)/2]}
    \gamma_{\,n,\,k} \,(xy)^{k} (x+y)^{n-2k+1}.
\end{equation}

The coefficients $\gamma_{\,n,\,k}$ are called
the $\gamma$-coefficients of the Eulerian polynomials.
Remarkably, Foata and Sch\"utzenberger discovered
a combinatorial interpretation of the coefficients
$\gamma_{\,n,\,k}$, that is,
for $n\geq 1$ and $1 \leq k \leq [(n+1)/2]$,
$\gamma_{\,n,\,k}$ equals
the number of permutations of $[n]$ with $k$ descents, but no double descents. Here a double descent of
a permutation $\sigma\in S_n$ is defined to be
an index $1 \leq i \leq n-1$ such that $\sigma_{i} > \sigma_{i+1} > \sigma_{i+2}$.

The nonnegativity of the coefficients
$\gamma_{\,n,\,k}$ has been referred to as the $\gamma$-positivity. This property of the Eulerian polynomials
and other polynomials along with the $q$-analogues has been extensively
 studied ever since,
see, for example,~\cite{
Athanasiadis-2018,
Branden-2008,
Chow-2008,
Foata-Schutzenberger-1970,
Foata-Strehl-1974,
Gasharov-1998,
Han-Ma-2020,
Lin-2017,
Lin-Zeng-2015,
Ma-Ma-Yeh-2019,
Ma-Ma-Yeh-2019-B,
Ma-Ma-Yeh-Yeh-2021,
Ma-Yeh-2017,
Petersen-2015,
Shapiro-Woan-Geta-1983,
Shin-Zeng-2012,
Stembridge-1992,
Yan-Zhou-Lin-2019}.

A context-free grammar is a set of substitution rules on a set of variables $X$.
A variable can be substituted with a polynomial (or a Laurent polynomials) in $X$.
Our starting point is the grammar of Dumont for the
Eulerian polynomials, namely,
\[ G= \{ x\rightarrow xy, \quad y \rightarrow xy\},\]
which can be expressed as a differential operator
\begin{equation} \label{Delta}
\Delta = xy {\partial \over \partial x} + xy
               {\partial \over \partial y}.
\end{equation}
Ma-Ma-Yeh~\cite{Ma-Ma-Yeh-2019} realized that
by a change of variables
\[ u=xy, \quad v=x+y,\]
the above grammar is transformed into a new grammar
\begin{equation}
H=\{ u \rightarrow uv, \quad v \rightarrow 2u\},
\end{equation}
ensuring the $\gamma$-positivity of the
Eulerian polynomials.

Without the setting of a grammar, the above argument can
be recast in terms of the differential operator $\Delta$
 in (\ref{Delta}). Clearly, we have
 $\Delta(x+y) = 2xy$ and $\Delta(xy) =xy (x+y)$. Since
 $\Delta$ is a derivative, we
 see that
for $n\geq 1$, $A_n(x,y)=\Delta^{n-1}(xy)$ is
 a polynomial in $x+y$ and $xy$ with nonnegative coefficients.

It turns out that the grammar $H$ plays an essential role in
the combinatorial understanding of the $\gamma$-coefficients
of the Eulerian polynomials.
First, we notice that the
 grammar $H$ bears a striking resemblance  to the following grammar for
0-1-2 increasing trees and the Andr\'e polynomials, namely,
\[ G=\{ x \rightarrow xy, \quad y \rightarrow x\}. \]

Recall that for $n\geq 1$, a 0-1-2 increasing tree on
$[n]$ is a
rooted increasing tree on $[n]$ for which
every vertex has at most two children. For $n\geq 1$, the
  Andr\'e polynomials are defined by
\[ E_n(x,y) = \sum_{T} x^{l(T)} y^{u(T)},\]
where the sum ranges over 0-1-2 increasing trees $T$
on $[n]$, $l(T)$ denotes the number of
leaves of $T$ and $u(T)$ denotes the number of vertices of $T$
having an only child.

Examining the factor of two in the grammar $H$, we are
guided precisely to the structure of 0-1-2 increasing
plane trees.
 This formulation is in agreement with
the known interpretation in reference to binary increasing trees
 on $[n]$
with exactly $k$ leaves and no vertices with left children
only. It is also in agreement with the formula of Han-Ma~\cite{Han-Ma-2020} in terms of
0-1-2 increasing trees on $[n]$ with  $k$ leaves.
Nevertheless, it seems to be convenient to work with 0-1-2 increasing
plane trees in order to describe the labeling consistent
with the grammar $H$.

The grammatical approach
associated with a grammatical
 labeling of 0-1-2 increasing plane trees offers a test ground for the main
result of this paper, which is concerned with the
trivariate second-order Eulerian polynomials on
Stirling permutations, introduced by Gessel and Stanley~\cite{Gessel-Stanley-1978}, see also Elizalde \cite{Elizalde-2021}.
For $n\geq 1$,  let $[n]_2$ denote the multiset
$\{1,1, 2, 2, \ldots, n, n\}$.
  A permutation $\sigma=\sigma_1\sigma_2 \cdots
\sigma_{2n}$ of $[n]_2$ is
said to be a Stirling permutation if for any $1\leq j \leq n$ the elements
between the two $j$'s in $\sigma$, if any, are greater than $j$. For $n\geq 1$, the set of Stirling permutations
of $[n]_2$ is denoted by $Q_n$.
A descent and an ascent of $\sigma \in Q_n$ can be defined
analogously to the case of an ordinary permutation.
For a Stirling permutation $\sigma$, we
adopt the convention that $\sigma$ is patched with
a zero both at the beginning and at the end, that is,
$\sigma_0=\sigma_{2n+1} = 0$.
The number of Stirling permutations of $[n]_2$
with $k$ descents is called the second-order
Stirling number, denoted by $C(n,k)$, or
$\left\langle \! \left\langle {n \atop k} \right\rangle
 \!\right\rangle$.

B\'ona~\cite{Bona-2009} introduced the notion of
a plateau of $\sigma \in Q_n$, which is
defined to be a pair of two adjacent elements
$(\sigma_{i}, \sigma_{i+1})$ such that
$\sigma_{i}=  \sigma_{i+1}$. More precisely,
 for
$\sigma \in Q_n$, the number plateaux, denoted
 ${\rm plat}(\sigma)$, is defined to be the number
 of indices $1\leq i \leq 2n$ such that
 $\sigma_{i } = \sigma_{i+1 }$.
B\'ona showed that for $n\geq 1$,
the statistics ${\rm asc}(\sigma)$, ${\rm des}(\sigma)$
and ${\rm plat}(\sigma)$ have the same distribution
over $Q_n$.
 Janson~\cite{Janson-2008} constructed an urn model
 to prove the symmetry of the joint distribution of
 the three statistics.

 It should be noted that a plateau of a
Stirling permutation was defined earlier by Dumont~\cite{Dumont-1980}
in the name of a repetition.
For $n \geq 1$, Dumont defined the
 polynomials $C_n(x,y,z)$ as
\[ C_n(x,y,z) = \sum_{\sigma\in Q_n} x ^{{\rm des } (\sigma)}
y ^{{\rm asc } (\sigma)}  z ^{{\rm plat } (\sigma)}.  \]
Note that for $n \geq 1$ and any $\sigma \in Q_n$,
\begin{equation} {\rm des } (\sigma) +
{\rm asc } (\sigma) +
{\rm plat } (\sigma)= 2n+1.
\end{equation}

 Dumont~\cite{Dumont-1980} obtained the following recurrence relation.
For $n\geq 1$,
\begin{equation} \label{r-c-n}
C_{n+1}(x,y,z)=xyz \left( {\partial \over \partial x}
+ {\partial \over \partial y}
+ {\partial \over \partial z} \right) C_n(x,y,z)
\end{equation}
with $C_1(x,y,z)=xyz$.
A refinement of the above recurrence
relation (\ref{r-c-n}) was  established     by Haglund-Visontai~\cite{Haglund-Visontai-2012}.

The differential operator in the recurrence relation
(\ref{r-c-n})
can be prescribed as a grammar
\begin{equation} \label{Grammar-HV}
 G=\{ x \rightarrow xyz, \quad y \rightarrow xyz, \quad
         z\rightarrow xyz\}.
         \end{equation}
Indeed, the grammatical labeling
as described in the arXiv version of~\cite{Chen-Hao-Yang-2021} or in~\cite{Ma-Ma-Yeh-2019} is essentially the same
argument as that given in~\cite{Haglund-Visontai-2012},
which in turn is in the same vein as the
recursive construction of  Janson~\cite{Janson-2008}.
So the above grammar $G$ in~(\ref{Grammar-HV})
should be
attributed to Dumont~\cite{Dumont-1980}.

The symmetry of $C_n(x,y,z)$ suggests that we may
 consider the expansion into the elementary
symmetric functions, as denoted by
\[u= x+y+z, \quad v = xy+xz+yz, \quad w =xyz.\]

 Thanks to the idea
 of Ma-Ma-Yeh~\cite{Ma-Ma-Yeh-2019}, we come to the following
grammar
\begin{equation}
H=\{u \rightarrow 3w, \quad v \rightarrow 2uw, \quad
      w \rightarrow vw  \}.
\end{equation}
Utilizing this grammar, we realize that the
argument for the $\gamma$-expansion of the Eulerian
polynomials can be carried over to the expansion of
$C_n(x,y,z)$ into the elementary symmetric functions.
To be more specific, we prove that for $n\geq 1$,
 $C_n(x, y, z)$ is a polynomial in $u,v, w$ whose
  coefficients can be interpreted in terms of
  0-1-2-3 increasing plane trees.

Very recently, Ma-Ma-Yeh-Yeh~\cite{Ma-Ma-Yeh-Yeh-2021-B}
have extended the grammatical approach to the $e$-positivity of the
symmetric
polynomials in $k+1$ variables defined on $k$-Stirling
permutations. For $k\geq 1$, they defined the polynomials $C_n(x_1, x_2, \ldots, x_{k+1})$ by means of the grammar $G=\{x_i \rightarrow x_1x_2\cdots
x_{k+1} \, |  \, i =1,2, \ldots, k+1\}$ and provided
 a grammatical labeling for $k$-Stirling permutations,
 which makes the symmetry property transparent.
It turns out that the polynomials $C_n(x_1,x_2, \ldots, x_{k+1})$
 are associated with the generating function of the joint
 distribution of the numbers of ascents,  $j$-plateaux and descents
 of $k$-Stirling permutations of $\{1^k, 2^k, \ldots, n^k\}$, introduced by Janson-Kuba-Panholzer~\cite{JKP-2011}, where $i^k$ represents
 $k$ occurrences of $i$. The symmetry
 of $C_n(x_1,x_2, \ldots, x_{k+1})$ was discovered in~\cite{JKP-2011}.
 The
 coefficients in the $e$-expansion of $C_n(x_1, x_2, \ldots, x_{k+1})$
  can be interpreted in terms
 of
increasing plane trees for which no vertex
has more than $k+1$ children, see~\cite{Ma-Ma-Yeh-Yeh-2021-B}.

  The background on the use of context-free grammars
  for combinatorial enumeration including the
  notion of a grammatical labeling can be found in~\cite{Chen-1993, Chen-Fu-2017}.
  In the next section, we shall give a glimpse of how
  to compute a generating function based on a context-free
  grammar. In a certain sense, this approach can be thought
  of as a formal calculus in the spirit of
   the symbolic method, while we may
  enjoy  the advantage that
   there is no fear of the lack of rigor.

\section{A grammatical calculus for $A_n(x,y)$}

A context-free grammar
can also be understood as a formal differential operator. For the purpose of combinatorial enumeration, the
variables are attached to combinatorial structures,
whereas the rules reflect the recursive construction
of combinatorial objects. Computationally speaking,
a grammar is  a derivative which is often informative
for deriving the generating functions.

Let us take the Eulerian polynomials $A_n(x,y)$ to
demonstrate the efficiency of the grammatical calculus.
Dumont~\cite{Dumont-1996}
 discovered the following grammar for  $A_n(x,y)$:
\begin{equation} \label{gxy}
G=\{x\rightarrow xy, \quad y\rightarrow xy.\}
\end{equation}
Let $D$ denote the formal derivative with respect to the
above grammar $G$. Dumont
 showed that $A_n(x,y)$
can be generated by the grammar $G$, that is, for $n\geq 1$,
\begin{equation}\label{a-d}
A_n(x,y) = D^n(x).
\end{equation}
Chen and Fu~\cite{Chen-Fu-2017}
 introduced the notion of a
grammatical labeling in the sense that
the grammar $G$
preserves information of significance along with the generation of
permutations in $S_{n+1}$ from permutations
in $S_{n}$.

If we express the formal derivative in terms of a
differential operator, the above relation (\ref{a-d}) can be written as
$A_1(x,y)= xy$ and for $n\geq 1$,
\[ A_{n+1}(x,y) =  xy \left( {\partial \over \partial x}
                      + {\partial \over
                       \partial y}  \right)
                       A_{n}(x,y).
\]
The above recurrence relation also appeared in Haglund-Visontai~\cite{Haglund-Visontai-2012}.
It is apparent that
for $n\geq 1$,
\begin{equation}
 A_n(x) = A_n(x,y)|_{y=1}.
\end{equation}

Let us proceed to present a derivation of the
well-known generating function of $A_n(x)$:
\begin{equation}
\label{Gen-A-n-x}
\sum_{n\geq 0} A_n(x) {t^n \over n!} = {1-x\over 1-xe^{(1-x)t}}.
\end{equation}

As noted by Carlitz and Scoville~\cite{Carlitz-Scoville-1974}, it is not so easy to
 recover the above generating function from the
 recurrence relation for the Eulerian numbers.
 However, by employing
  the grammatical calculus
one can perform this task with ease.
 In fact, we find it more convenient to
deal with the following generating function of
the bivariate version $A_n(x,y)$.

\begin{thm} \label{thm-a-n}
 Set $A_0(x,y)=y$. Then we have
\begin{equation} \label{g-a-x-y}
\sum_{n=0}^\infty A_n(x,y){t^n \over n!} =
{y-x \over 1 - xy^{-1} e^{(y-x)t} }.
\end{equation}
\end{thm}

It is evident that setting $y=1$ in
(\ref{g-a-x-y}) yields (\ref{Gen-A-n-x}).
To present a grammatical proof of (\ref{g-a-x-y}),
recall that for a Laurent polynomial $f$ in $x$ and $y$, the generating function of $f$ with respect to the grammar
$G$ is defined by
\begin{equation} \label{Gen-d}
{\rm Gen}(f, t) = \sum_{n\geq 0} D^n(f) {t^n \over n!}.
\end{equation}
Assume that $g$ is also a Laurent polynomial in $x$ and $y$. The first and foremost property of $D$ is that
it is a derivative, that is,
\begin{equation}
D(fg) = D(f) g + f D(g),
\end{equation}
and hence it obeys the Leibniz rule
\begin{equation}
D^n(fg) =  \sum_{k=0}^n {n \choose k} D^{k}(f) D^{n-k}(g),
\end{equation}
for any $n\geq 0$. This implies the multiplicative property
\begin{equation} \label{Gen-fg}
{\rm Gen}(fg, t) = {\rm Gen}(f, t)\, {\rm Gen}(g, t) .
\end{equation}

\noindent
{\it Proof of Theorem \ref{thm-a-n} by Using the
Grammar of Dumont.} Under the assumption
$A_0(x,y)=y$, we have $A_n(x, y)=D^n(y)$.
So our goal is to compute the
generating function ${\rm Gen}(y, t)$.

For the formal derivative $D$ with respect to the
 grammar $G$ in (\ref{gxy}), we have
\begin{equation} \label{D-x-1}
D(y^{-1} ) =  -  y^{-2} D(y) = -x y^{-1}
\end{equation}
and
\begin{equation} \label{D-x-1-A}
 D(xy^{-1})=xy^{-1} (y-x).
 \end{equation}
As noted in~\cite{Chen-Fu-2017}, since $x-y$
is a constant with respect to $D$, we deduce that
for $n\geq 0$,
\begin{equation} \label{D-x-1-y}
D^n(xy^{-1} ) = x  y^{-1} (y-x)^n.
\end{equation}
In light of the property (\ref{Gen-fg}) and the fact $D(c)=0$ when $c$ is a constant, it suffices
to consider ${\rm Gen}(y^{-1},t)$, since
\begin{equation} \label{G-y-2}
{\rm Gen}(y, t) = {1 \over {\rm Gen}(y^{-1},t)}.
\end{equation}
Using (\ref{D-x-1}), we obtain that
\begin{equation}
{\rm Gen}(y^{-1},t) = \sum_{n\geq 0}
D^n(y^{-1}){t^n \over n!}
= y^{-1} - \sum_{n\geq 1} D^{n-1} (xy^{-1}) {t^n \over n!}.
\end{equation}
Invoking (\ref{D-x-1-y}), we deduce that
\begin{eqnarray*}
{\rm Gen} (y^{-1},t) & = & y^{-1} -
\sum_{n\geq 1}  xy^{-1}(y-x)^{n-1}{t^n \over n!} \\[9pt]
& = & y^{-1} -
   {xy^{-1} \over y-x} \left( e^{(y-x)t} -1 \right) \\[9pt]
& = & {1 - xy^{-1} e^{(y-x)t} \over y-x},
\end{eqnarray*}
which completes the proof by utilizing (\ref{G-y-2}).   \qed

A variation of the generating function of $A_n(x,y)$
was considered by Carlitz-Scoville~\cite{Carlitz-Scoville-1974},
which equals
\begin{equation}
F(t)= { e^{xt} - e^{yt}\over  xe^{yt} -y e^{xt}}.
\end{equation}
Recall that $A_0(x,y)$ is
defined to be $y$ and $A_n(x,y)$ is symmetric in
$x$ and $y$ for $n\geq 1$. We are led to consider
the generating function of $A_n(x,y)$ for $n\geq 1$.
Using the above generating function,
we see that
\begin{equation}
 {\rm Gen} (y, t) -y = xy \,{ e^{xt} - e^{yt}\over  xe^{yt} -y e^{xt}}.
\end{equation}
The reason for the appearance
of the factor $xy$ in the above
expression becomes evident once we
take a close look at the definition
of $F(t)$ given by Carlitz-Scoville~\cite{Carlitz-Scoville-1974}.

\section{The  $\gamma$-positivity of $A_n(x,y)$}

The $\gamma$-coefficients $\gamma_{\,n,\,k}$ of the
Eulerian polynomials $A_n(x,y)$ as given in
(\ref{A-n-gamma}) and (\ref{gamma-x-y})  have a number of combinatorial interpretations.
For the purpose of this paper, we shall single out the one
in  connection with 0-1-2 increasing plane trees.

A 0-1-2 increasing plane tree on $[n]$ is an
increasing plane tree for which each vertex has
degree at most two, where the degree of a vertex is
referred to the number its children.
For a 0-1-2 increasing plane
tree $T$ on $[n]$, assume that it has $f_0$ leaves and
 $f_2$vertices of degree two, then it is easily seen that
\begin{equation}\label{f20}
 f_2=f_0-1.
 \end{equation}
 Let $s(n,k)$ be the number of
0-1-2 increasing trees on $[n]$ with $k$ leaves,
and let $t(n,k)$ be the number of 0-1-2 increasing
plane trees on $[n]$ with $k$ leaves. Then we have
\begin{equation} \label{st}
t(n,k) = 2^{k-1} s(n,k).
\end{equation}

We now turn to the observation of Ma-Ma-Yeh~\cite{Ma-Ma-Yeh-2019} on a grammatical
explanation of the $\gamma$-positivity of $A_n(x,y)$.
Observe that
\[
D(xy) =  (x+y)xy, \quad
D(x+y)  =   2xy.\]
If we set
\[  u=xy, \quad v=x+y, \]
then we get $D(u)=uv$ and $D(v)=2u$.
In other words, we have a new grammar
\begin{equation}\label{GH}
 H= \{ u \rightarrow uv, \quad v \rightarrow 2u\}.
\end{equation}
Let $D$ denote the formal derivative with respect to the
grammar $G$ as well as the grammar $H$. It is safe to do so
since $G$ and $H$ have distinct variables.
Since for $n\geq 1$,
\[ A_n(x,y)=D^n(x) =D^{n-1}(u),\]
we infer that $A_n(x,y)$ is a polynomial
in $u$ and $v$ with nonnegative coefficients.
That is to say, the polynomials $A_n(x,y)$
are $\gamma$-positive.

Let $T$ be a 0-1-2 increasing plane tree on $[n]$,
where $n \geq 1$. We define a grammatical labeling of $T$ as follows.
A leaf is labeled by $u$, a degree one vertex is labeled
by $v$ and a degree two vertex is labeled by $1$.
 The weight of $T$ is defined to be the product of the
 labels associated with the vertices of $T$.

 For example, Figure 1 is a 0-1-2 increasing plane on $[6]$ with weight $u^3v$, where the grammatical labels are in
 parentheses.
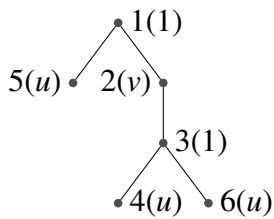
\begin{figure}[!ht]
\begin{center}
\begin{tikzpicture}
\node [tn,label=0:{$1$}($1$)]{}[grow=down]
	[sibling distance=12mm,level distance=8mm]
  child {node [tn,label=180:{$5$}($u$)]{}
            edge from parent
                node[left=2pt ]{}
                }
        child {node [tn,label=180:{$2$}($v$)](three){}
            [sibling distance=12mm,level distance=8mm]
	       child{node [tn,label=0:{$3$}($1$)](two){}
            child{node [tn,label=0:{$4$}($u$)]{}}
            child{node [tn,label=0:{$6$}($u$)]{}
                edge from parent
                    node[right=2pt ]{}
            }}
            edge from parent
                node[left=2pt ]{}
            }

        ;
\end{tikzpicture}
\caption{A 0-1-2 increasing plane tree.}
\end{center}
\end{figure}

Increasing plane trees are also called plane recursive trees, see Janson~\cite{Janson-2008}.
The above grammatical labeling of 0-1-2
increasing plane trees shows that the
$\gamma$-coefficients for the Eulerian polynomials
can be interpreted based on 0-1-2 increasing plane
trees.

\begin{thm} For $n\geq 1$ and $1\leq k \leq [(n+1)/2]$,
the number $\gamma_{\,n,\,k}$
equals the number of 0-1-2 increasing plane trees
on $[n]$ with $k$ leaves.
\end{thm}

It is not hard to transform a 0-1-2 increasing
plane tree into
 a permutation without double descents.
There is a one-to-one correspondence $\phi$ between the
set of permutations on $[n]$ with
$k$ descents and no double descents and the
set of 0-1-2 increasing plane trees on
$[n]$ with $k$ leaves.
This is exactly the
classical bijection between permutations and
 increasing binary trees, restricted to binary
  trees without vertices
 having only left children,  see Stanley~\cite{Stanley-I}.

The structure of 0-1-2 increasing plane trees can
be employed to partition the set of permutations of $[n]$ into classes
similar to the classification according to the
Foata-Strehl action on permutations defined based on the notion of $x$-factorization~\cite{Foata-Strehl-1974},
see also Br\"and\'en~\cite{Branden-2008}. Foata and Strehl gave a proof of \eqref{A-n-gamma} using the group action on $S_n$, by which $S_n$ can be partitioned into equivalence classes.

 Let $T$ be a 0-1-2 increasing plane tree on $[n]$.
We now consider a labeling of $T$ by assigning
a label $x$ or $y$ to a degree one vertex, and a label
$xy$ to a leaf, and $1$ to a degree two vertex. Let
$\alpha(T)$ denote the set of labeled trees obtained
from $T$, and let $w(T)$ denote the sum of weights
of all trees in $\alpha(T)$. Assume that
$T$ has $k$ leaves. As  has already been mentioned, $T$ has
$k-1$ degree two vertices. Thus it contains
$n+1-2k$ degree one vertices. It follows that the
total weight of the trees in $\alpha(T)$ amounts
to
\begin{equation} \label{wt}
w(T) = (xy)^k (x+y)^{n+1-2k}.
\end{equation}

The above relation reveals
 that the set of permutations
of $[n]$ can be partitioned into classes with each class
corresponding to a 0-1-2 increasing plane tree $T$.
More precisely, for a 0-1-2 increasing plane tree $T$,
we get a set $\alpha(T)$
  of labeled trees with each corresponding to a permutation. Now, the sum of weights of these
   labeled trees is given by (\ref{wt}). It is
   readily seen that the sum of weights of trees
   can be  translated into a weighted sum of permutation statistics.

We notice that
 0-1-2 increasing plane trees with the
above labeling scheme can be represented
as increasing binary trees.
For each labeled tree in $\alpha(T)$, we can
represent it by an increasing binary tree on $[n]$.
For a degree one vertex $v$, if it is labeled by $x$,
then we turn its child into a left child as in a binary
tree, otherwise, we turn it into a right child as in a
binary tree. By the classical bijection between
permutations and increasing binary trees, we find that the labeling
of a 0-1-2 increasing plane tree is suitable for
keeping track of the number of descents of a permutation.
As can be seen,  the grammatical labeling
of 0-1-2 increasing plane trees provides a combinatorial
justification of the relation of Foata-Sch\"utzenberger
back to  the original form.

\section{The  Second-order Eulerian Polynomials}

Gessel and Stanley~\cite{Gessel-Stanley-1978}
introduced the notion of Stirling permutations and
 defined the second-order
Eulerian polynomials $C_n(x)$ by $C_0(x)=1$
and for $n\geq 1$,
\begin{equation}
C_n(x)= \sum_{k=1}^{n} C(n,k)x^k,
\end{equation}
where $C(n,k)$ is the number of Stirling permutations
on $[n]_2$ with $k$ descents.
A homogeneous version of $C_n(x)$ is given by
\begin{equation}
C_n(x,y) = \sum_{k=1}^{n} C(n,k)x^k y^{2n+1-k}.
\end{equation}

Let $C_n(x,y,z)$ be the trivariate polynomials
first defined by Dumont~\cite{Dumont-1980} and rediscovered by Janson~\cite{Janson-2008}.
As a symmetric function in $x,y,z$, $C_n(x,y,z)$ can be
expressed as a polynomial in the elementary
symmetric functions in $x,y,z$. If the coefficients
 are all nonnegative, we say that
the symmetric function is $e$-positive, see Stanley~\cite{Stanley-II}. We shall show that for $n\geq 1$,
 $C_n(x,y,z)$ is $e$-positive along with a
 combinatorial interpretation of the coefficients.

 Let $G$ be the following grammar
\begin{equation}\label{G-xyz}
G=\{ x\rightarrow xyz, \quad y\rightarrow xyz, \quad
     z\rightarrow xyz\}.
\end{equation}
Let $D$ denote the formal derivative with respect
to $G$. It has been shown by Dumont~\cite{Dumont-1980} that
for $n\geq 1$,
\begin{equation}
  C_n(x,y,z) = D^n(x).
\end{equation}

For $n\geq 1$, assume that
\begin{equation}\label{C-n-g}
C_n(x,y,z) =\sum_{i+2j+3k=2n +1} \gamma_{\,n,\,i,\,j,\,k} (x+y+z)^i
(xy+xz+yz)^j (xyz)^k .
\end{equation}

Let \[
u=x+y+z,\quad  v=xy+xz+yz, \quad  w=xyz. \]
Then we have
\begin{equation} \label{dh}
 D(u) = 3w , \quad
 D(v) = 2uw, \quad
 D(w) = vw.
 \end{equation}

 For $n\geq 1$, we have $C_n(x,y,z)=D^{n-1}(xyz)$,
  where $D$ is the formal derivative with respect to the
  grammar $G$ in (\ref{G-xyz}). Now, we might as well use the same symbol $D$ for the formal derivative
  with respect to the grammar $H$.
   For $n\geq 1$, we have \[ D^n(x)=D^{n-1}(w),\]
   which is clearly a polynomial in $u,v,w$ with
   nonnegative coefficients.  In other words, $C_n(x, y, z)$ is $e$-positive.

The main objective of this paper is to give a combinatorial interpretation
 of the coefficients $\gamma_{\,i,\, j,\, k}$ in
 (\ref{C-n-g}). The relations in (\ref{dh}) prompt us to define the grammar
 \begin{equation} \label{Grammar-uvw}
 H=\{ u \rightarrow 3w, \quad v \rightarrow 2uw, \quad
      w \rightarrow vw\}.
 \end{equation}

 \begin{thm} For $n \geq 1$ and $i+2j+3k=2n+1$,
 the coefficient $\gamma_{\,i,\,j,\, k}$ in the expansion
  (\ref{C-n-g}) of $C_n(x, y, z)$ equals the
 number of 0-1-2-3 increasing plane trees on $[n]$
 with $k$ leaves, $j$ degree one vertices and $i$
 degree two vertices.
 \end{thm}

\proof
Let $T$ be a 0-1-2-3 increasing plane tree on
$[n]$. We first give a labeling of $T$   as
follows. Label a leaf by $w$,
a degree one vertex by $v$,  a degree two vertex by $u$ and a degree three vertex by $1$. Given any 0-1-2-3 increasing plane tree $T$ on $[n]$ with $k$ leaves, $j$ degree one vertices and $i$ degree two
vertices, it has $n-i-j-k$ vertices of degree three. Taking
the number of edges into consideration, we get
\[
 3(n-i-j-k) + 2i + j = n-1.
\]
Thus we have verified that \[ i+2j+3k=2n+1.\]

Let us examine
how to generate a 0-1-2-3 increasing
plane tree $T'$ on $[n+1]$ by adding $n+1$ to
 $T$ as a leaf. We can   add $n+1$ to $T$
 only  as a child of a vertex $r$ that is not of degree three.
 Thus we have the following three possibilities.

\noindent Case 1:
The vertex $r$ is a leaf with label $w$.
 In the resulting tree $T'$,
$r$ becomes a degree one vertex with label $v$ and $n+1$
 becomes a leaf with label $w$. This operation
 corresponds to the substitution $w \rightarrow vw$.

\noindent Case 2:
The vertex $r$ is a degree one vertex with label $v$.
In this case, $n+1$ can be attached to $r$ either
as the first child, or  the second child. In either case,
in the resulting tree $T'$,
$r$ becomes a degree two vertex with label $u$ and $n+1$
 becomes a leaf with label $w$. This operation
 corresponds to the substitution $v \rightarrow 2uw$.

 \noindent Case 3:
The vertex $r$ is a degree two vertex with label $u$.
In this case, $n+1$ can be attached to $r$ either
as the first child, or the second child, or the third
 child. In either case,
in the resulting tree $T'$,
$r$ becomes a degree three vertex with label $1$ and $n+1$
 becomes a leaf with label $w$. This operation
 corresponds to the substitution $u \rightarrow 3w$.

The aforementioned three cases exhaust all the possibilities to
construct a 0-1-2-3 increasing plane tree $T'$ on $[n+1]$ from
a 0-1-2-3 increasing plane tree $T$ on $[n]$ by
adding $n+1$ as a leaf. Since each case corresponds to an
application of a substitution rule in $H$, we see that
for $n\geq 1$,
$D^n(x)$ equals the sum of the weights of
0-1-2-3 increasing plane trees on $[n]$, that is,
\begin{equation}
D^n(w)= \sum_{i+2j+3k =2n+1} \gamma_{\,i,\,j,\, k}\, u^i v^j w^k.
\end{equation}
Therefore, $\gamma_{\,i,\,j,\,k}$ equals the number
of 0-1-2-3 increasing plane trees on $[n]$ with
$k$ leaves, $j$ degree one vertices and $i$ degree two
vertices. \qed

Figure 2 is an illustration of a 0-1-2-3 increasing
plane tree on $[10]$, where the grammatical labels are in
  parentheses.

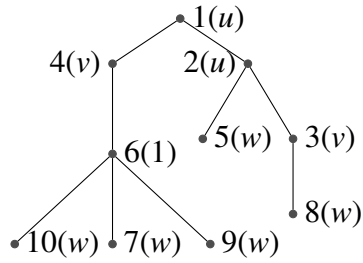
\begin{figure}[!ht]
\begin{center}
\begin{tikzpicture}
\node [tn,label=0:{$1$}($u$)]{}[grow=down]
	[sibling distance=18mm,level distance=6mm]
child {node [tn,label=180:{$4$}($v$)](four){}
            [sibling distance=13mm,level distance=12mm]
	       child{node [tn,label=0:{$6$}($1$)](three){}
            child{node [tn,label=0:{$10$}($w$)]{}}
            child{node [tn,label=0:{$7$}($w$)]{}}
            child{node [tn,label=0:{$9$}($w$)]{}
                edge from parent
                    node[right=2pt ]{}
            }}
            edge from parent
                node[left=2pt ]{}
            }
             child {node [tn,label=180:{$2$}($u$)]{}
                [sibling distance=12mm,level distance=10mm]
                child{node [tn,label=0:{$5$}($w$)]{}}
                child{node [tn,label=0:{$3$}($v$)]{}
                    child{node [tn,label=0:{$8$}($w$)]{}
                        edge from parent
                              node[right=2pt ]{}
                              }
                   edge from parent
                         node[right=2pt ]{}
                         }
            edge from parent
                node[left=2pt ]{}
                }

        ;
\end{tikzpicture}
\caption{A 0-1-2-3 increasing plane tree.}
\end{center}
\end{figure}

In view of  the above
combinatorial interpretation or the relation $D^n(u)=D(D^{n-1}(u))$, we conclude with the
 following recurrence relation:
\begin{equation}
\gamma_{\,n,\,i,\,j,\,k}=
3(i+1)\gamma_{\,n-1,\,i+1,\,j,\,k-1}
+2(j+1)\gamma_{\,n-1,\,i-1,\,j+1,\,k-1}
+k\gamma_{\,n-1,\,i,\,j-1,\,k}
\end{equation}
with $\gamma_{\,1,\,0,\,0,\,1}=1$ and $\gamma_{\,1,\,i, \,j, \,k}=0$ if $k \neq 1$.

\vskip 5mm \noindent{\large\bf Acknowledgments.}  We wish to thank M. B\'ona, D. Foata, S.-M. Ma, P. Paule and the referee for invaluable comments and suggestions.
This work was supported
by the National
Science Foundation and the Ministry of Science and
 Technology of China.

\end{document}